\documentclass[11pt,reqno]{article}

\usepackage{amsmath,amsthm,amsfonts,amssymb}
\usepackage{t1enc}
\usepackage[ps2pdf]{hyperref}
\usepackage{graphicx}
\usepackage{epsfig}
\usepackage{graphpap}
\let\phi=\varphi

\chardef\bslash=`\\ 

\newtheorem{thm}{Theorem}[section]
\newtheorem{conj}{Conjecture}[section]

\newtheorem{lem}[thm]{Lemma}
\newtheorem{prop}[thm]{Proposition}

\theoremstyle{definition}
\newtheorem{defn}{Definition}[section]

\theoremstyle{remark}

\newtheorem{que}{Question}[section]
\newtheorem{ex}{Example}[section]

\def\orig{o}

\newcommand{\x}{\mathbf{x}}
\newcommand{\y}{\mathbf{y}}

\newcommand{\Z}{\mathbb{Z}}
\newcommand{\N}{\mathbb{N}}
\renewcommand{\P}{\mathbb{P}}
\newcommand{\E}{\mathbb{E}}
\newcommand{\T}{\mathbb{T}}

\newcommand{\eval}[2][\right]{\relax
  \ifx#1\right\relax \left.\fi#2#1\rvert}

\title{A note on spider walks}
\author{Christophe Gallesco$^{1}$ \and Sebastian M\"uller$^{2}$ \and
Serguei Popov$^{3}$}

\begin{document}
\maketitle {\footnotesize \noindent $^{~1}$ Instituto de Matem{\'a}tica e Estat{\'\i}stica, Universidade de
S{\~a}o Paulo, rua do Mat{\~a}o 1010, CEP 05508--090, S{\~a}o Paulo, SP, Brazil

\noindent e-mail: \texttt{gallesco@ime.usp.br}

\noindent $^{~2}$  Institut f\"{u}r Mathematische Strukturtheorie, Technische Universit\"{a}t Graz, Steyrergasse
30, 8010 Graz, Austria

\noindent e-mail: \texttt{dr.sebastian.mueller@gmail.com}

\noindent $^{~3}$  Department of Statistics, Institute of Mathematics, Statistics and Scientific Computation, University of Campinas--UNICAMP, rua S\'ergio Buarque de Holanda 651, CEP 13083--859, Campinas SP, Brazil

\noindent e-mail: \texttt{popov@ime.unicamp.br}, \quad \noindent url:
\texttt{http://www.ime.unicamp.br/$\sim$popov}

}

{\abstract Spider walks are systems of interacting particles. The particles move independently as long as their movements do not violate some given rules describing the relative position of
the particles; moves that violate the rules are not realized. The goal of this paper is to study qualitative properties, as recurrence, transience, ergodicity, and positive rate of escape of these Markov processes. }

\section{Introduction}

Let us start with an informal description of a particular spider walk. Imagine there are two particles performing  simple random walks on a graph in continuous time.
These particles are tied together with a rope of a certain length, say $s$. As long as the rope is not tight their movements are independent. If the rope is tight (particles have a distance $s$ of each other) the rope prevents the particles to jump \emph{away} from each other.

More generally we can think of a spider walk as a system of~$k$, $k\in \N$, interacting particles that move
independently according to some Markov process as long as their movement does not violate some given rules
concerning their relative positions. If a move of a particle violates a rule the particle stays at its position
and waits until the particle jumps to another location or the movement of another particle change the relative
position of the particles. In other words, for each particle in some position $x$ and each neighbouring site $y$
there is an exponential clock with rate $q(x,y)$ independent of the rest of the process. Say the clock
associated to a particle in $x$ and to the edge $xy$ rings then the particle moves to $y$ if the new position
accords with the rules and stays in $x$ otherwise. A formal construction is given in
Subsection \ref{subsec:spidergraph}.

This note addresses questions about the qualitative characteristics of the original random walk and the spider
walk. In particular, we are interested if the \emph{rope} may change properties like recurrence vs.\ transience,
ergodicity vs.\ non-ergodicity and positive speed vs.\ zero speed. The first observation that one may make is
that a spider walk can be described as a random walk of one particle on some appropriate graph that we call the
spider graph. Some of the above questions may then be answered in comparing the two graphs. For example consider
the simple random walk (SRW) as the underlying Markov process. Then, if the original graph and the spider graph
are roughly isometric, the SRW is recurrent if and only if the spider walk is recurrent. In the more general
setting of reversible Markov processes one can compare the two processes by dint of rough embeddings, see
Definition \ref{defn:embedding} and Theorem \ref{th:re_trans}. Let us denote $Q^S$ for the Markov process on the
spider graph $G^S$, i.e., for the interpretation of the spider walk as a \emph{one particle walk}. If $Q$ and
$Q^S$ are roughly equivalent the original process and the spider walk are either both recurrent or both
transient. This is the case if the underlying process and the \emph{rules} are transitive, see Theorem
\ref{thm:tr_re}. The same holds true if we leave the setting of transitivity but assume that the conductances
are bounded, see Theorem \ref{thm:re_trans_2}. If we drop the hypotheses on transitivity and bounded
conductances a transient Markov process can bear a recurrent spider, see Example \ref{ex:tree_rooted_at_inf},
and a transient spider can originate from a recurrent Markov process, see Example \ref{ex:re-tr}.

There is no  analogue to Theorem \ref{thm:tr_re} treating positive- and null-recurrence since there exists no
positive recurrent quasi-transitive Markov chain. (This follows from Theorem~1.18, Lemma~3.25 and Theorem~3.26
in \cite{woess}). In other words, every transitive spider is null recurrent or transient. In the Examples
\ref{ex:pos-null}, \ref{ex:null-pos}, and \ref{ex:tree_rooted_at_inf} we describe situations where the random
walk and the spider walk have different ergodic behaviours.

A natural follow up question is whether the rate of escape is positive or zero. For this let us first consider SRW on graphs with a positive anchored isoperimetric constant. It is known that positivity of the anchored isoperimetric constant implies positivity of the speed. By observing that the anchored isoperimetric constant of $G$ is positive if and only if the one of the spider graph is positive we obtain for this class of graphs that the speed of the spider is positive. We believe that on transitive graphs the speed of the SRW is positive if and only if the one of the spider is positive. This is not true for transitive random walks in general: see Example~\ref{ex:tree_rooted_at_inf} where the underlying random walk has positive speed
but the spider walk has zero speed and \emph{vice versa}. Furthermore, we study spiders with bounded span on the integers (with drift) and on homogeneous trees and observe two different qualitative behaviours: while the speed
$V(s)$ of the spider  on the line converges to the speed of the random walk as $s$ goes to infinity, the speed
of the spider walk on the tree converges to $0$ as the span increases. In both cases the speed of the spider
walk is strictly smaller than the speed of the underlying process. While this is not true in general, see Example
\ref{ex:tree_rooted_at_inf}, we conjecture it to hold for SRW on transitive graphs. We conclude with some
questions concerning structural properties of graphs and the behaviour of the spider walk.

Our model and results have also a motivation coming from evolutionary dynamics and molecular cybernetics of
multi-pedal molecular spiders, compare with \cite{antal:07}. There, different mappings between various models of
spiders and simple excursion processes are established. We also want to mention \cite{GMPV:09} where the spider
walk in random environment on~$\Z$ is studied.

Let us comment on some nomenclature. A Markov chain (process) is called a random walk if the process is \emph{somehow} adapted to the structure of the graph. Furthermore, we switch freely between the continuous time and discrete time version (jump chain) of the process. Due to the description of our particle system as \emph{spider}, the particles are frequently called \emph{legs}. In most of the basic definitions and results concerning random walks and reversible Markov chains we follow the two monographs \cite{lyons:book} and \cite{woess} where more details and references can be found.

\section{Spider walk}
\subsection{Definition}\label{subsec:spidergraph}
Let $G=(V,E)$ be a rooted, undirected, connected graph with vertex set $V,$ edge set $E$, and root $\orig$.  As usual $d(\cdot,\cdot)$ denotes the graph distance and $\deg(x)$ the degree of the vertex $x$. We often identify the graph with its vertex set and write $x\in G$ for a vertex $x$. A Markov process $(X_t)_{t\geq 0}$ starting in $\orig$ is defined on $G=(V,E)$ through a transition matrix $Q=(q(x,y))_{x,y\in G}.$  We always assume that the rates are bounded and that $G$ and $Q$ are adapted, i.e., $q(x,y)=q(y,x)=0$ if and only if $xy\not\in E$. This assumption in particular implies, since the graph is connected, that the Markov process is irreducible.

We define the spider walk with $k$ legs in a very general setting. We distinguish between the different legs and describe the spider walk through $S_{t}=(S_{1,t},\ldots,S_{k,t})$ where $S_{i,t}$ stands for the position of the
$i$th leg at time $t$. For each $x\in V$ let us fix $L(x)$ which is a finite subset of
$\Lambda(x):=\{(x_1,\ldots, x_k):~x_1=x,~x_i\in V\}$. We call $L(x)$ the set of local configurations of the
spider at position $x$ and write $\x=(x_1,\ldots, x_k)$ for its elements.  We define the transition rates $Q^S$
from $\x=(x_1,\ldots,x_k)$ to $\y=(y_1,\ldots,y_k)$ as follows:
\begin{equation*}
q^S(\x,\y)=\left\{ \begin{array}{cl}
  q(x_i, y_i) & \mbox{if there exists exactly one index}~i~\mbox{such that}~x_i\neq y_i  \\
  0 & \mbox{otherwise.} \\
\end{array}\right.
\end{equation*}

Together with some initial position $\bar \x$ of the spider, the sets $L(x),~x\in V$, and the transition matrix $Q$ define the spider
walk $(S_{t})_{t\geq 0}$ through $Q^S$. Furthermore, we denote $G^S=(V^S, E^S)$, where
$$ V^S=\bigcup_x L(x)\mbox{ and } E^S:=\{(\x,\y):~q^S(\x,\y)>0\}.$$ We call $G^S$ the spider graph of the spider walk $(S_{t})_{t\geq 0}$. Its root is
$\bar\x$. Sometimes it will be convenient to write the vertices of the spider graph also as $\{ (\ell_{i}(x)):~i\in\{1,\ldots, |L(x)|\},~x\in G\},$ where $1,\ldots, |L(x)|$ corresponds to some enumeration of the set $L(x)$.

Since in this general definition the spider walk is not necessarily irreducible (even though the process $Q$ is irreducible) we will concentrate on two types of spider walks: namely spider with bounded span and transitive spiders.

\subsection{Spider with bounded span} We consider a spider walk with
$k$ legs and span bounded by $s$, i.e., the maximal distance between the $k$ legs does not exceed $s$. For each $x$  let
$$L(x)=\{(x_1,\ldots,x_k):~x_1=x,~ \max_{i,j} d(x_i,x_j)\leq s,~ x_i\neq x_j\}.$$ As a starting
position $\bar \x$ we may choose $\bar x_1,\ldots, \bar x_k$ as a non-intersecting path of length $k$ starting
from $\orig$. Observe hereby that we do not allow two legs to be at the same position.

\begin{ex}\label{ex:1}
Simple random walk on $\Z$ and $2$-leg spider with $s=3.$ We assume the first leg to be the leftmost leg. In
this case $L(x)=\{(x,x+1), (x,x+2), (x,x+3)\}$, compare with Figure \ref{fig:ex:1}. A part of the corresponding
spider graph is drawn in Figure \ref{fig:spidergraph}.
\begin{figure}[h]
\begin{center}
\includegraphics[scale=0.8]{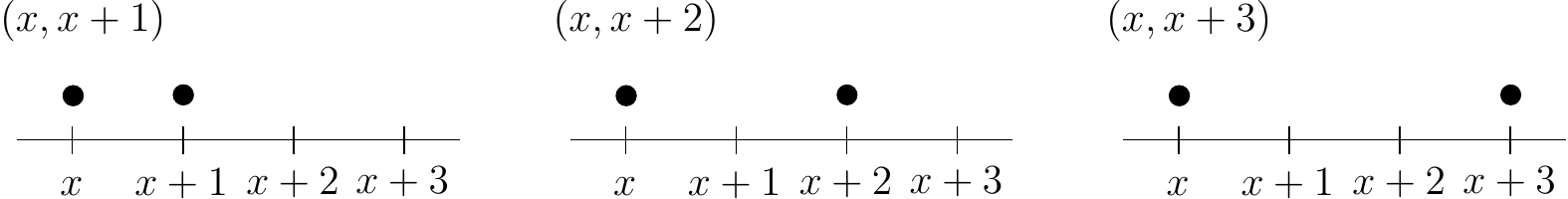}
\caption{The three possible local configurations for a spider with bounded span ($k=2,~s=3$) at position
$x$.}\label{fig:ex:1}
\end{center}
\end{figure}

\begin{figure}[h]
\begin{center}
\includegraphics[scale=0.8]{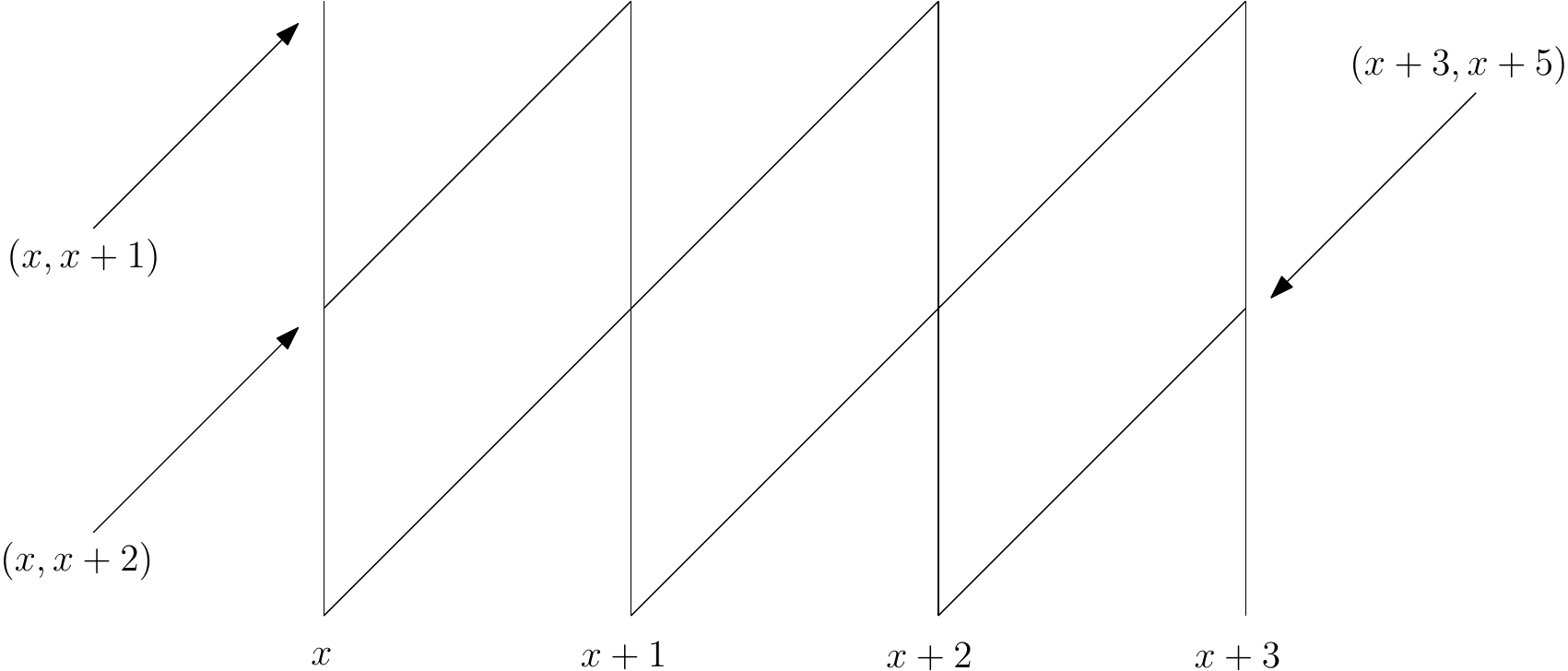}
\end{center}\caption{Part of the spider graph of a spider with bounded span ($k=2,~s=3$)}\label{fig:spidergraph}
\end{figure}
\end{ex}

Note that a spider walk  with bounded span is \emph{a priori} not irreducible; e.g. the spider on $\Z$ with $2$ legs and span bounded by $1$.  A natural assumption that ensures irreducibility is $k\leq s$.

\subsection{Transitive spider}
Let us recall the definition of transitive graphs. A graph $G=(V,E)$ is said to be transitive if the
automorphism group $AUT(G)$ acts transitively on $G$, i.e., for all $x,y\in G$ there exist a $\gamma\in AUT(G)$ such that $y=\gamma x$. Informally spoken: the graph $G$ looks the same from every vertex. Let $Q$ be a transition matrix describing a Markov process on $G$ and $AUT(G,Q)$ be the group of all $\gamma\in AUT(G)$ which
satisfy $q(\gamma x, \gamma y)=q(x,y)$ for all vertices $x,y$ in $G$. We say the Markov process $Q$ is transitive if the group $AUT(G,Q)$ acts transitively on $G$.  
We can extend this idea to the set of local configurations $L(x)$.
Define $AUT(G,Q,L)$ as the group of all $\gamma\in AUT(G,Q)$ which satisfy $L(\gamma x)=L(x)$ for all $x\in G$. We say the spider $(G,Q,L)$ is transitive if  $AUT(G,Q,L)$ acts transitively on $G$.

\begin{ex}\label{ex:transitivespider1}
Let us define a spider with $3$ legs on the line $\Z$ with the following set of local configurations:
\[ L(x)=\{ (x, x+1,x+2), (x, x+1, x+3), (x, x+2, x+3), (x,x+2,x+4)  \}.\] This defines a transitive spider since $L(x)$ is just the translation of $L(0)$. Observe, that this spider is not of bounded span since we excluded the local configurations $(x,x+1,x+4)$ and $(x,x+3,x+4)$.
\end{ex}

\begin{ex}\label{ex:transitivespider2}
Every spider with bounded span on a transitive graph is a transitive spider. For example consider the $2$-leg spider with bounded span $s=2$ on the direct product $G$ of $\Z_{3}$ and $\Z$. Elements of $G$ are written as a tuple $(u,x)$ with $u\in \Z_{3}$ and $x\in \Z$. The set of local configurations can then be written as
\begin{eqnarray*}
L((u,x))& = & \{ ((u,x),(u\pm 1,x)), ((u,x),(u\pm 2,x)),\cr & &~~((u,x),(u,x\pm 1)), ((u,x),(u-1,x\pm 1)), ((u,x),(u+1,x\pm 1))  \}.   
\end{eqnarray*} In Figure \ref{fig:transitive:spider} the position of the first leg is labelled by the black ball and the possible positions of the second leg are indicated by the grey balls.

\begin{figure}[h]
\begin{center}
\includegraphics[scale=0.8]{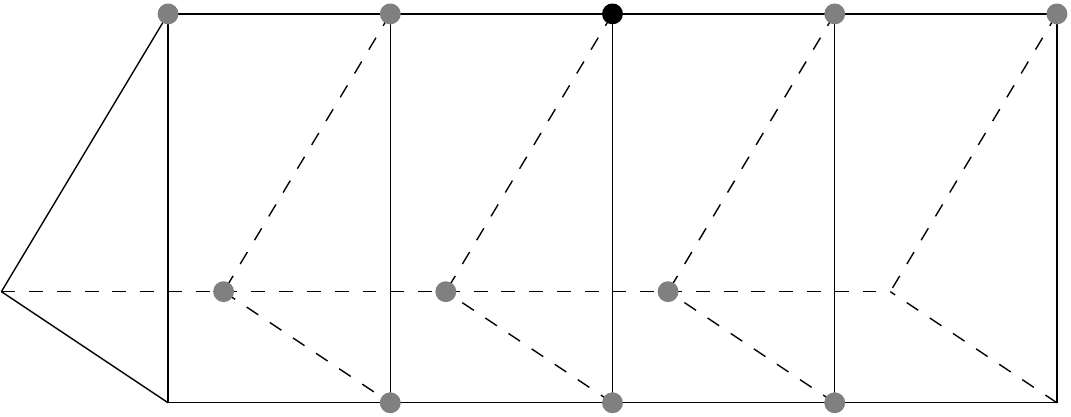}
\end{center}\caption{Set of local configurations of  a  spider with bounded span ($k=s=2$) on the direct product  $\Z_{3}\times \Z$.}\label{fig:transitive:spider}
\end{figure}

\end{ex}

 Let $\{\ell_1,\ldots, \ell_{|L(\orig)|}\}$ be an arbitrary enumeration of $L(\orig)$. 
 This induces an enumeration on the local configuration at position $x$ as follows: choose some $\gamma\in AUT(G,Q,L)$ with $\gamma \orig =x$ and enumerate $L(x)$ such that $\gamma \ell_{i}(\orig)=\ell_{i}(x)$.

A graph is called quasi-transitive if the automorphism group acts with finitely many orbits. Recall when a group $\Gamma$ acts on a set $G$, the group orbit of an element $x\in G$ is defined as
$\Gamma x:=\{\gamma x:~ \gamma\in\Gamma\}$. Let $O_i$ be the orbits of $AUT(G,Q,L)$ on $G$. The vertex set of
the factor graph consists of the orbits and two orbits $O_i$ and $O_j$ are connected by an edge if there exists
$u\in O_i$ and $v\in O_j$ such that $uv\in E$. We can define the factor chain on the factor graph by $\tilde
p(O_i,O_j)=\sum_{w\in O_j} p(x,\omega)$, where $x\in O_i$ is arbitrary.

Due to the above definition, the spider graph $G^S$ of a transitive spider is a quasi-transitive graph and hence
its geometry can be compared with those of the underlying graph $G$. One possibility of comparing two different
graphs is the concept of rough-isometry. Equip two graphs $G=(V,E)$ and $G'=(V',E')$ with their natural metrics
$d$ and $d'$. A mapping $\phi: V \to V'$ is called rough isometry if there are positive constants $\alpha$ and
$\beta$ such that
\begin{equation*}\label{eq:roughiso}
\alpha^{-1} d(x,y) - \beta \leq d'(\phi(x), \phi(y))\leq \alpha d(x,y) +\beta\quad \forall x,y\in
V\end{equation*} and such that every vertex in $G'$ is within distance $\beta$ from the image of $G$. In this
case we say that the two graphs (or metric spaces) $(G,d)$ and $(G',d')$ are roughly isometric.

\begin{lem}\label{lem:trans_ri}
Let  $(G,Q,L)$ be a transitive spider walk. Then $G$ and the spider graph $G^S$ are roughly isometric.
\end{lem}

\begin{proof}
Denote $d^S$ the usual graph distance in $G^S.$ Define $\phi: V \to V^S,~x \mapsto \ell_1(x).$ Clearly,
$d(x,y)\leq d^S(\ell_1(x), \ell_1(y))$ for all $x,y\in V$. Due to the transitiveness we have for all $x,y\in V$
that $d^S(\ell_1(x), \ell_1(y))\le\alpha d(x,y)$ with $\alpha=d^S(\ell_1(v), \ell_1(w)),$ where $vw\in E.$
Furthermore, there is some constant $\beta$ such that $d^S(\ell_1(x),\ell_i(x))\leq \beta$ for all $1\leq i\leq
u$ and all $x\in G.$
\end{proof}

Lemma \ref{lem:trans_ri} does not hold for general graphs, compare with the following Example
\ref{ex:no:quasi-iso}.

\begin{ex}\label{ex:no:quasi-iso}
Consider the $1$-dimensional grid $\Z$ and attach to each vertex of the form $2^k, k\in\N,$ an additional
vertex. The spider graph of the $2$-leg spider with span $s=2$ is not quasi-isometric to the underlying graph
since $d^S(\ell_i(x),\ell_j(x))$ is not bounded for some $\ell_i, \ell_j$ and $x$. To see this let
$x=2^k+2^{k-1}$ and $\ell_i$ be a configuration  where the first leg is to the left of the second and $\ell_j$
where it is to the right of the second leg.
\end{ex}

\section{Recurrence and transience}
In order to define recurrence and transience of a Markov process $(X_t)_{t\geq 0}$ (on a discrete state space
$X$) it is convenient to pass to its jump chain $(Y_n)_{n\geq 0}$. The jump chain is the (discrete) sequence of states visited by the continuous-time Markov process $(X_t)_{t\geq 0}$.  Its transition probabilities $P=(p(x,y))_{x,y\in X}$  are $p(x,y):=- q(x,y)/ q(x,x)$ for $x\neq
y$ and $0$ if $x=y$. We also can write this relation as $Q=R(P-I)$ where $R$ is the diagonal matrix whose diagonal
entries are $\{q(x,x):~x\in X\}$ and $I$ the diagonal matrix with diagonal entries equal to $1$.

An irreducible Markov chain $(Y_n)_{n\geq 0}$ is called recurrent if $\P(\exists n\geq 1:~Y_n=y\mid Y_0=x)=1$
for all $x,y\in X$, otherwise it is called transient. We say the Markov process is recurrent (respectively,
transient) if its jump chain is recurrent (respectively, transient). We say a Markov process is reversible if
there exists a positive vector $\mu$ such that the detailed balance condition, $\mu(x) q(x,y)=\mu(y) q(y,x)$,
holds for all $x,y\in X$. The reversibility of the Markov process carries over to its jump process: let
$\pi(x)=-q(x,x)\mu(x)$, then  $\pi$ is the reversible measure of $(Y_n)_{n\geq 0}$, i.e.,  $\pi(x)
p(x,y) = \pi(y) p(y,x)$ for all $x,y\in X$. In what follows we restrict us to the study of reversible Markov
processes. Notice that in general (if the underlying process is not reversible) the spider walk may develop
\emph{singular} behaviour, e.g., the first leg may return infinitely many times to the starting position while
the spider walk itself is transient. Furthermore, the interpretation of a Markov chain as an electrical network
is restricted to reversible Markov chains: any reversible Markov chain defines an electric network with
conductances $c(x,y):=\pi(x) p(x,y)$. Note that due to the reversibility we have $c(x,y)=c(y,x)$ and hence the conductances
$(c(x,y))_{xy\in E}$ can be seen as weights on the edges of the graph $G$ according to which a random walker
chooses its next position. The resistance of an edge is defined as $r(x,y):=1/c(x,y)$.

Let us first turn to transitive spiders.  As an immediate consequence of Lemma \ref{lem:trans_ri} we have that
if the simple random walk (SRW), i.e., $p(x,y)=1/\deg(x)$ if $xy\in E$ or equivalently $q(x,y)=1 $ if $xy\in E,$
on a transitive graph $G=(V,E)$ is transient then every transitive spider on this graph is transient and if some
spider is transient then the SRW is transient as well. This fact follows from the well-known result that rough
isometries preserve recurrence and transience, compare with Theorem 3.10 in \cite{woess}. In order to generalize
this result to transitive spiders we use the concept of rough embeddings.

\begin{defn}\label{defn:embedding}
Let $G_1$ and $G_2$ be electrical networks with resistances $r^{(1)}$ and $r^{(2)}$. We say that a map $\phi$
from the vertices of $G_1$ to the vertices of $G_2$ is a rough embedding if there are constants $a,b< \infty$
and a map $\Phi$ defined on the edges of $G_1$  such that
\begin{enumerate}
 \item[a)] for every edge $xy$ in $G_1,$ $\Phi(xy)$ is a non-empty simple
oriented path of edges in $G_2$ joining $\phi(x)$ and $\phi(y)$ with
$$ \sum_{e\in \Phi(xy)} r^{(2)}(e)\leq a r^{(1)}(x,y),$$
where the sum is over all edges in the path $\Phi(xy).$

\item[b)] $\Phi(xy)$ is the reverse of $\Phi(yx)$

\item[c)] for every edge $e$ in $G_2$ there are no more than $b$ edges in $G_1$ whose image under $\Phi$
contains~$e.$
\end{enumerate}

We call two networks roughly equivalent if there are rough embeddings in both directions.

\end{defn}

There is the result of \cite{kanai:86} stating that the type is preserved under rough embeddings.

\begin{thm}\label{th:re_trans}
If there is a rough embedding from $G_1$ to $G_2$ and $G_1$ is transient, then $G_2$ is transient.
\end{thm}

We are now able to prove the following

\begin{thm}\label{thm:tr_re}
Let $Q$ be a reversible and transitive Markov process on $G$ and $(G,Q,L)$ be an irreducible transitive spider.
Then, $Q$ is transient if and only if $(G,Q,L)$ is transient.
\end{thm}
\begin{proof}
First we recall the following general fact. Let $(X,P)$ be a transitive reversible Markov chain with reversible
measure $\pi.$ For $x,y\in X$ let $n$ such that $p^{(n)}(x,y)>0$. Then $\pi(x) p^{(n)}(x,y)=\pi(y)p^{(n)}(y,x)$
and $\pi(\gamma x) p^{(n)}(x,y)=\pi(\gamma y)p^{(n)}(y,x)$, with $\gamma\in AUT(X,P).$  Dividing yields, $\frac{
\pi(\gamma x)}{\pi(x)}=\frac{\pi(\gamma y)}{\pi(y)}$ and hence the function $g(\gamma):=\frac{\pi(\gamma
x)}{\pi(x)}$ does not depend on $x.$  Consequently, $g$ is an exponential on $AUT(X,P),$ i.e., $g(\beta
\gamma)=g(\beta) g(\gamma).$ Moreover, the function
$f(\gamma):=\frac{r(\gamma(x),\gamma(y))}{r(x,y)}=\frac1{g(\gamma)}$ is an exponential.

Let us use the above observation for our setting. Notice that both the electrical network of the underlying Markov chain $P$ and the one of the spider walk $P^S$ are restrictions of the electrical network $P^{I}$ describing the movement of $k$ independent particles. 
Denote by $G^{I}$ the graph corresponding to $P^{I}$.
Since $G$ is transitive, $G^{I}$ is transitive with corresponding automorphism group $AUT(G^{I})$ (take automorphisms coordinatewise translation). We can compare the resistances in $G$ and $G^S$ as induced subgraphs of $G^{I}$. Clearly, if some edge $e$ is in $G$
and $G^S$ then $r(e)=r^S(e)=r^{I}(e).$ 
Let $\ell_l(x) \ell_m(y) $ be an edge  in $G^S$. Then there exist some $\gamma\in  AUT(G^{S})$, $v\sim \orig$ and $i,j$ such that $\ell_{l}(x)=\gamma \ell_{i}(\orig)$ and $\ell_{m}(y)=\gamma \ell_{j}(v)$. Hence,

\begin{equation}\label{eq:res_comp}
 \frac{r^{S}(\ell_l(x),\ell_m(y))} {r(x,y)}= \frac{r^{S}(\gamma(\ell_i(\orig)),\gamma(\ell_{j}(v)))} {r(\gamma(\orig),\gamma(v))}=\frac{f(\gamma)r^S(\ell_i(\orig),\ell_j(v))}{f(\gamma) r(\orig,v)}:=k(i,j,v).
\end{equation}

We are now ready to apply Theorem \ref{th:re_trans}. To do this we construct a rough embedding $\Phi$ from
$G_1=G$ to $G_2=G^S.$ Let $\phi: V \to V^S$ be $\phi(x):=\ell_1(x).$ In order to construct $\Phi$ we fix some
reference points $v, w \in V$ with $vw\in E$ and let $\Phi(vw)$ be some (arbitrary but fixed) shortest path from
$\ell_1(v)$ to $\ell_1(w)$ in $G^S.$ For $xy\in E$ we define $\Phi(xy)$ as the shortest path from $\ell_1(x)$ to
$\ell_1(y)$ such that $\gamma(\Phi(xy))=\Phi(vw)$ for some  $\gamma\in AUT(G^S)$. Due to this construction we have that $\Phi(xy)$ is the reverse of $\Phi(yx)$ and due to the quasi-transitivity of $G^S$ we
have that for every edge in the factor graph there are only finitely many edges in $G$ whose image under $\Phi$
contains~$e$. We have to check the first property:
$$ \sum_{e\in \Phi(xy)} r^{(S)}(e)\leq a r(x,y),$$ for some
$a>0.$ But this holds with $a:=|\Phi(xy)|\cdot \max k(i,j,v)$ using (\ref{eq:res_comp}).

It remains to construct a rough embedding from $G^S$ to $G.$ Let $\phi(l_i(x)):=x$ and $\Phi(l_i(x)l_j(y)):=xy$, then one may verify the three properties of Definition \ref{defn:embedding} as in the first part.
\end{proof}

We now turn to spider walks with bounded span where we use the following fact that is left as an exercise in \cite{lyons:book}, Proposition 2.17.

\begin{prop}\label{prop:ri_re}
Let $G_1$ and $G_2$ be two infinite roughly isometric networks with conductances $c_1$ and $c_2.$ If $c_1, c_2,
c_1^{-1}, c_2^{-1}$ are all bounded and the degrees in $G_1$ and $G_2 $ are all bounded, then $G_1$ is roughly
equivalent to $G_2.$
\end{prop}

\begin{thm}\label{thm:re_trans_2}
Let $Q$ be a reversible  Markov process with bounded conductances and $Q^S$ an irreducible $k$-leg spider walk
with span $s$. Then, the Markov process is transient if and only if the spider is transient.
\end{thm}
\begin{proof}
Due to Example \ref{ex:no:quasi-iso} we can not use the spider graph in order to show that there is a
quasi-isometry between the two processes. We need a different encoding of the position of the spider walk. To do
this we choose an enumeration of the vertices $G$ in such a way that the root $\orig$ of $G$ corresponds to $0.$
The position of the spider is now defined as the closest (in graph distance) position of a leg to the origin. If there are several closest positions  we choose the one with the smallest number in the enumeration. Analogously to Subsection \ref{subsec:spidergraph} we can define another \emph{spider graph} with global positions $x\in G$ and the set of local configurations that we again denote by $L(x)=\{\ell_i(x):1\leq
i\leq |L(x|)\}$. The fact that the conductances of the network of the spider walk are bounded follows from the
fact that it is a subnetwork of the network describing the movement of $k$ independent particles, compare with
the proof of Theorem \ref{thm:tr_re}. Due to Proposition \ref{prop:ri_re} it remains to show that $G$ and the
new spider graph are roughly isometric.

First,  we show that the distance between  two local configurations $\ell_i(x)$ and $\ell_j(x)$ of the same
global position is uniformly bounded. Recall that a local configuration in $x$ can be described as the
sequence $x=x_1,\ldots,x_k$. We call a configuration \emph{lined} if $d(x_i,x_{i+1})=1$ for $1\leq i < k$.
Observe that it takes at most $k(s+k)$ steps to get from $\ell_i(x)$ to any lined configuration in $x$.
Consequently, since this procedure is invertible we obtain that
\[ d^S(\ell_i(x), \ell_j(x))\leq 2k(s+k),\] where $d^S$ is the graph distance in the new spider graph.
 Second, observe that for each $x$ and $y$ there exists some $i^*$ and $j^*$ such that
  $d^S(\ell_{i^*}(x),\ell_{j^*}(y))\leq k d(x,y)$ and
 hence
 \[ d^S(\ell_i(x), \ell_j(y))\leq k d(x,y) + 4k(s+k),~\forall \ell_i(x), \ell_j(y).\]
Eventually, rough isometry follows as in the proof of Lemma \ref{lem:trans_ri} .
\end{proof}

If we drop the hypotheses on transitivity in Theorem~\ref{thm:tr_re} or bounded conductances in Theorem~\ref{thm:re_trans_2} a transient Markov chain can bear a recurrent spider, compare with Example
\ref{ex:tree_rooted_at_inf} in Section~\ref{sec:speed}, and a transient spider can originate from a recurrent
Markov chain, compare with Example~\ref{ex:re-tr}.

\begin{ex}\label{ex:re-tr}
\emph{Recurrent Markov chain and transient spider walk.} \newline\noindent We consider an example of a Lamperti random walk, that is a nearest neighbour random walk $(Y_{n})_{n\geq 0}$  on $\N$ with  asymptotic zero drift. The mean drift of $(Y_{n})_{n\geq 0}$ is defined as $ \mu(x)  = E [ Y_{n+1} - Y_{n} \mid Y_n=x ]$ and is supposed to go to $0$ as $n$ goes to $\infty$. There is the following criterion, due to \cite{lamperti:60},  for recurrence and transience  in terms of the mean drift, see Theorem~3.6.1 (i)-(ii) in~\cite{fayolle:95} : If there exists a number $B$ such that $\mu(x)\leq \frac{1}{2x}$ for $x\geq B$ then the Markov chain $(Y_{n})_{n\geq 0}$ is recurrent. On the other hand if for some $B$ and $\theta>1$ we have $\mu(x)\geq \frac{\theta}{2x}$ for $x\geq B$, then the Markov chain is transient.

Since $(Y_{n})_{n\geq 0}$  is a nearest neighbour walk we have $q(x,x-1)=1-q(x,x+1)\in(0,1)$ and hence that $\mu(x)=2 q(x,x+1)-1$. Letting $q(x,x+1)=(2x+1)/(4x)$ it follows that the Markov chain is recurrent since $\mu(x)=1/(2x)$.

Now, consider the $2$-leg spider with span $s=2$. We assume the first leg to be the left leg of the spider. In
this case $L(x)=\{(x,x+1), (x,x+2)\}$, see Figure~\ref{fig:ex:re-tr}.
\begin{figure}[ht]
\begin{center}
\includegraphics[scale=0.8]{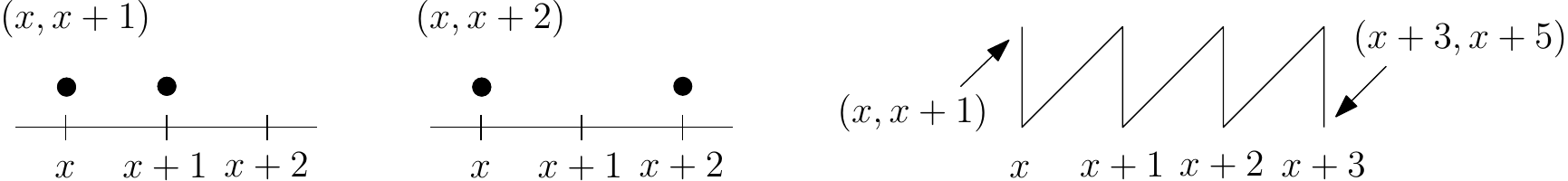}
\caption{The two possible local configurations for a spider walk with bounded span ($k=2,~s=2$) at position $x$ and a
part of the spider graph.}\label{fig:ex:re-tr}
\end{center}
\end{figure}

Let us identify $(x,x+i)$ with $2x+i-1.$ Hence the spider graph can be seen as a \emph{stretched line}, compare
with Figure \ref{fig:ex:re-tr}, and consequently the spider walk is itself a nearest neighbor random walk on the
line. We calculate its mean drift as
$$\mu^s(2x)= \frac{q(x+1, x+2)-q(x, x-1)}{q(x+1, x+2)+q(x, x-1)}=\frac{2x+1}{4x(x+1)-1}$$ and
$$\mu^s(2x+1)= \frac{q(x, x+1)-q(x+2, x+1)}{q(x, x+1)+q(x+2, x+1)}=\frac{x+1}{2x(x+2)+1}.$$
Eventually, we obtain $\mu^s(x)\sim 1/x$ which implies transience of the spider.
\end{ex}

The two following examples demonstrate that for reversible Markov chains that are not quasi-transitive both can happen: positive recurrent Markov chain and null recurrent spider, and null recurrent Markov chain and positive recurrent spider.

\begin{ex}\label{ex:pos-null}
\emph{Positive recurrent Markov chain and null recurrent spider walk.}
\newline\noindent Let $G_N,$ $N\in\N$, be the line segment $[0,1,\ldots,N]$. The graph $G$ is constructed in identifying the $0$'s
of the graphs $G_N$. In order to distinguish the vertices of the different $G_N$ we denote $x_N$ for the
vertices of $G_N$. On $G$ we define the following Markov chain $(Y_n)_{n\geq 0}$:
$$ p(0,1_N)=\left(\frac12\right)^N,\quad
p(N_N,N_N)=1-p(N_N,(N-1)_N)=p,$$
$$p(x_N,(x+1)_N)=1-p(x_N,(x-1)_N)=p \quad \forall 1\leq x<N,$$
where $p\in(1/2,1)$. Let $T_0:=\min\{n>0:~Y_n=0\}$ be the first return time to $0$. It is straightforward to
show that
$$\E[T_0 \mid Y_0=1_N]\sim \left(\frac{p}{q}\right)^N,$$ with $q=1-p$. Observe that
the $2$-leg spider with span $s=2$ on $G_{N}$ behaves like a random walk with drift on the line segment $G_{2N}$, compare with Example \ref{ex:re-tr}.  Let $T^S_0$ be the first time the spider visits $(0,1_1)$. Hence
we obtain for the spider walk, denoted by $S_n$, that
\[\E[T^S_0 \mid S_0=(0,1_N)]\sim \left(\frac{p}{q}\right)^{2N}.\]
Now choosing $p/q=\sqrt{2}$ we obtain a positive recurrent Markov chain and a null recurrent spider.
\end{ex}

\begin{ex}\label{ex:null-pos}
\emph{Null recurrent Markov chain and positive recurrent spider walk.}
\newline\noindent We consider the nearest neighbor Lamperti random walk with asymptotic zero drift, compare with Example
\ref{ex:re-tr}. Recall the corresponding criterion for ergodicity, Theorem~3.6.1(iii)-(iv) in~\cite{fayolle:95}: If there
exists a number $B$ such that $0\geq \mu(x)\geq -\frac{1}{2x}$ for $x\geq B$ then the Markov chain $Y_n$ is null
recurrent. On the other hand if for some $B$ and $\theta>1$ we have $\mu(x)<- \frac{\theta}{2x}$ for $x\geq B$,
then the Markov chain is positive recurrent. Eventually, the Lamperti random walk with mean drift
$\mu(x)=-\frac1{2x}$ is null recurrent but the corresponding $2$-leg spider with span $s=2$ is positive
recurrent.
\end{ex}

\section{Speed}\label{sec:speed}

There is no analogue to Theorem \ref{thm:tr_re} treating positive and zero speed. Already for transitive spiders
it might be that the random walk has zero speed but the spider has positive speed and the random walk has
positive speed but the spider walk has zero speed, compare with Example~\ref{ex:tree_rooted_at_inf}. For the
latter phenomenon we also refer to \cite{GMPV:09} for an example in random environment.

\begin{ex}\label{ex:tree_rooted_at_inf}
We consider the homogeneous tree $\T=\T_M$ with degree $M$ and root $\orig$. A ray $\langle x_0, x_1,\ldots,
\rangle$ is an infinite path from $\orig$ to infinity that does not backtrack, i.e., $x_i\neq x_j$ for all
$i\neq j$. Two rays are said to be equivalent if their symmetric difference has finitely many vertices.  We call
the set of all equivalence classes of rays the (end) boundary of $\T$, denoted by $\partial \T$. If $x\in \T$
and $\xi\in \partial \T$ then $\xi$ has a unique representative which is a ray starting at $x,$ denoted by
$\langle x,\xi\rangle$. The confluent $x\curlywedge y$ of two vertices $x$ and $y$ with respect to a ray
$\omega$ is the first common vertex on the rays $\langle x,\omega\rangle$ and $\langle y,\omega\rangle$. For any
$x\in \T$ we define its height with respect to $\omega$ by
$$h(x)=d(x, x\curlywedge \orig)- d(\orig, x\curlywedge \orig),$$
where $d(\cdot,\cdot)$ is the natural graph distance. The $k$-th horocycle of $\T$ (with respect to $\omega$ and
$\orig$) is the set $H_k=\{x\in \T:~ h(x)=k\}$. For more details on this model we refer to Chapter 9  in
\cite{woess2} and to Figure \ref{fig:tree}. One can think of $\omega$ as the \emph{mythical ancestor} of the
\emph{genealogical} tree $\T$. Each $x\in H_k$ has exactly one neighbour (father) $x^-$ in $H_{k-1}$ and $M-1$
neighbours (sons), $y^-=x$, in $H_{k+1}$.

\begin{figure}
\setlength{\unitlength}{0.75cm}
\begin{picture}(10,10)
\put(1,1){\vector(1,1){10}} \put(10,10.5){\makebox(0,0)[b]{$\omega$}}

\put(1.5,1.5){\line(1,-1){0.5}}

\put(2.5,2.5){\line(1,-1){1.5}} \put(3.5,1.5){\line(1,-1){0.5}} \put(3.5,1.5){\line(-1,-1){0.5}}

\put(4.5,4.5){\line(1,-1){3.5}} \put(6.5,2.5){\line(-1,-1){1.5}} \put(5.5,1.5){\line(1,-1){0.5}}
\put(7.5,1.5){\line(-1,-1){0.5}} \put(4.4,4.6){\makebox(0,0)[b]{$\orig$}}

\put(8.5,8.5){\line(1,-1){3}}

\multiput(1.5,1.5)(0.25,0){55}{\circle*{0.005}} \put(0.5,1.5){\makebox(0,0)[l]{$H_2$}}

\multiput(1.5,2.5)(0.25,0){55}{\circle*{0.005}} \put(0.5,2.5){\makebox(0,0)[l]{$H_1$}}

\multiput(1.5,4.5)(0.25,0){55}{\circle*{0.005}} \put(0.5,4.5){\makebox(0,0)[l]{$H_0$}}

\multiput(1.5,8.5)(0.25,0){55}{\circle*{0.005}} \put(0.5,8.5){\makebox(0,0)[l]{$H_{-1}$}}
\end{picture}\caption{A part of the tree \emph{rooted} at $\omega$ and horocycles corresponding to $\orig$.}
\label{fig:tree}
\end{figure}
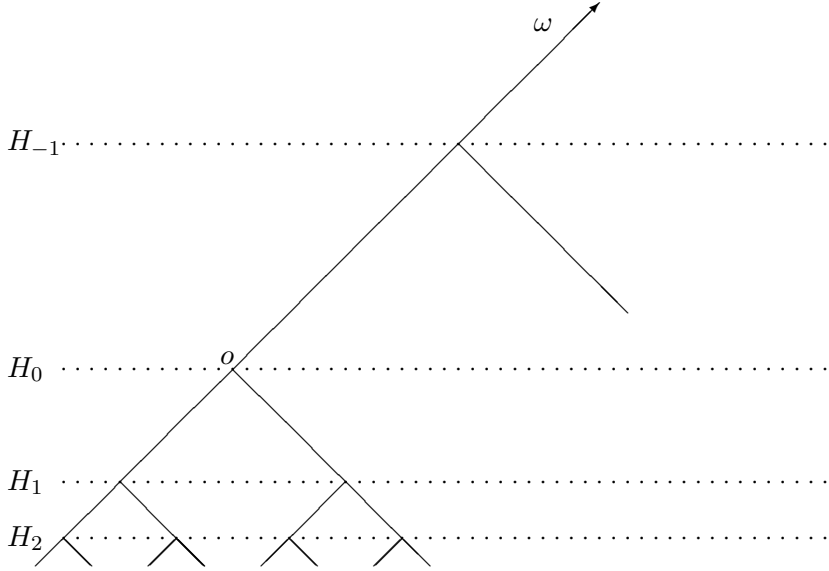

We define a nearest neighbour random walk  $(X_{n})_{n\geq 0}$ on $\T$:

$$p(x^-,x)=a/(M-1),\quad p(x,x^-)=1-a,$$ where $0<a<1$. It turns
out that the random walk is transient for all $a\in(0,1)$ and $\lim_{n\to\infty} |X_n|/n=|2a-1|$ a.s. The speed
of the jump chain $(S_n)_{n\geq 0}$ of the spider walk is calculated  using
\begin{equation}\label{eq:h(S_n)}
\lim_{n\to\infty} \frac1n h(S_n)=\lim_{n\to\infty} \frac1n\sum_{k=0}^{n-1} \left( H(S_{k+1})-H(S_k)\right).
\end{equation}

The averages in~(\ref{eq:h(S_n)}) are over an ergodic stationary sequence since the factor chain is ergodic and
hence the speed of the jump chain is
$$D=\sum_{i} \Pi(O_i) \E [ H(S_1)-H(S_0) \mid S_0\in O_i],$$ where
$\Pi(\cdot)$ denotes the stationary measure of the factor chain.

In order to obtain the speed in continuous time  we have to divide $D$ by the mean time it takes until a leg
jumps, i.e., by
$$T=\sum_{i} \Pi(O_i) \E [ \inf_{t>0}\{S_t\neq S_0\} \mid
\tau(S_0)=i].$$

Calculating the speed of the $2$-leg spider with $M=3$ and $a=1/2$ one observes: while for $s=1,2$ the speed
of the spider equals zero  the spider converges to $\omega$ with positive speed for $s=3$.  Now, consider
$\T_\orig$ as the subtree of $\T$ rooted in $\orig$ that consists of all vertices in $x\in H_k,~k\geq 0$ such
that $x$ lies on the geodesic $\langle \orig,\omega\rangle.$ In other words $\T_\orig$ is the genealogical tree
of $\orig$. We obtain for $s=3, M=4$, and $a=1/2$ that the random walk is null-recurrent on $\T_\orig$ while the
spider is positive recurrent.

Furthermore, it is possible to choose $a<1/2$ such that the spider has zero and the random walk positive speed.
In the same way, we can find some $a$ such that the random walk converges  to $\partial \T\setminus\{\omega\}$
with positive speed and the spider with positive speed to $\omega$. With this latter $a$ the random walk on
$\T_\orig$ is transient and the spider on $\T_\orig$ is positive recurrent.
\end{ex}

\subsection{Spider walk on the line}\label{subsec:line}
We consider the $2$-leg spider with span $s$ on $\Z$ with rates $q(x,x+1)=p,~q(x,x-1)=q,$ and $0$ otherwise. The
speed of the underlying Markov chain  clearly is $V=p-q.$ In order to compute the speed of the spider we first
calculate the stationary distribution $\Pi$ of the factor chain, compare with the calculation in Example
\ref{ex:tree_rooted_at_inf}:
$$\Pi(1)=\frac1{2s-2}, ~\Pi(i)=\frac1{s-1},~1<i<s,\mbox{ and } \Pi(s)=\frac1{2s-2}.$$ This yields to
$$ D(s)= \frac12 \frac{p-q}{p+q}\mbox{ and } T(s)= \frac{s}{2(s-1)} \frac1{p+q}.$$
Eventually, we obtain the speed of the spider
$$V(s)=\frac{D(s)}{T(s)} =(p-q) \left(1-\frac1s\right).$$
This  reproduces and generalizes  results in \cite{antal:07}, where the speed was calculated for
$p\in\{1/2,1\}$. Furthermore, notice that $V(s)<V$ for all $s\ge 2$ and $V(s)\to V$ as $s\to \infty$.

\subsection{Spiders walk on homogeneous trees}\label{subsec:tree}
In this subsection we consider SRW on homogeneous trees and transitive spiders. As mentioned in the introduction the speed of the spider will be positive if the anchored isoperimetric constant of the spider is positive. Furthermore, the speed of the spider is strictly smaller than the speed of the random walk:

\begin{thm}\label{thm:tree:slower_speed}
Let $Q$ be a  SRW on the homogeneous tree $\T_M$ of degree, $M\geq 3$, and speed $V$. Then any irreducible
transitive spider $(\T_M, Q, L)$ has positive speed strictly less than $V$.
\end{thm}
\begin{proof}
The idea of the proof is to compare the spider graph $G^S$ with a larger graph $H$ on which the SRW has the same
speed as the SRW on $G=\T_M$. Observe that in $G^S$ each vertex $\ell_i(x)$ has at most one neighbour in the set
$\{\ell_j(y)\}_{j\in L}$ for all $y\neq x$. We add edges to $G^S$ such that for all $x\in G, i\in L$ and all
$y\sim x$ there exists exactly one $j\in L$ such that $\ell_i(x)\sim \ell_j(y)$. In other words, in $H$ every
vertex $\ell_i(x)$ has exactly one neighbour in $\{\ell_j(y)\}_{j\in L}$ for all $y\neq x$. Clearly, the SRW (in
continuous time) on $H$ has the same speed as the SRW on $G$. It remains to show that the addition of edges
strictly increases the speed of the random walk. Let $\partial B_n$ be the set of vertices (of $G$) with distance $n$ to $\orig$.
Each vertex $x\in \partial B_n$ has one neighbour in $\partial B_{n-1}$ and $M-1$ neighbours in $\partial B_{n+1}$. Let us consider the set
of local configurations $\{\ell_i(x)\}_{i\in L}$ of the global position $x\in \partial B_n.$ For each edge that is added
to a local position in $x^-\in \partial B_{n-1}$ we add one corresponding edge to $y$ where $x\sim y$ and $y\in \partial B_{n+1}$.
Eventually, for each edge leading back to the origin $M-1>1$ edges are added that lead away from the origin. The
claim follows now by the observation that the stationary distribution of the factor chains (in continuous time)
of the local positions of $G_S$ and $H$ is the uniform distribution. The fact that the speed is positive follows
from the discussion in the introduction.
\end{proof}

Let us now consider the $2$-leg  spider with bounded span on the homogeneous tree of degree $d\geq 3.$ The
calculation of the speed becomes more complicated compared to those for the spider walk on $\Z$, see
Subsection \ref{subsec:line}. Observe that the number of local configurations is of order $s^2$ and exact
calculations of the stationary distribution of the factor chain become involved. Nevertheless, we are able to
prove the following asymptotic result that we conjecture to hold true for transitive spiders on homogeneous
trees in general.

\begin{thm}\label{conj:speed_tree}
Let $V(s)$ be the speed of a $2$-leg spider with span $s$ on a homogeneous tree $\T_{M}$ of degree $M\geq 3.$ Then,
$V(s)\to 0,$ if $s\to \infty.$
\end{thm}

\begin{proof}
For sake of simplicity we give the proof only for $M=3$. We encode the local configuration of the spider as the
tuple $(l,k)$, ($1\leq l,k\leq s$), where $l=d(x_1,x_2)$ is the distance  and $k=|d(x_1,\orig )-d(x_2,\orig)|$
is the height difference between the two legs $x_1$ and $x_2$. This factor chain has the following transition
rates $q$, compare with Figure \ref{figure1}.
\begin{itemize}
    \item $l,k=1$: $q((1,1),(2,2))=3$ and $q((1,1),(2,0))=1$
    \item $1<l=k<s$: $q((l,l),(l+1,l+1))=3,~q((l,l),(l+1,l-1))=1$
    and $q((l,l),(l-1,l-1))=2$
    \item $l,k=s$: $q((s,s),(s-1,s-1))=2$
    \item $0<l<s,k=0$: $q((l,0),(l-1,1))=2$ and
    $q((l,0),(l+1,1))=4$
    \item $l=0, k=s$: $q((0,s),(1,s-1))=2$
    \item $0<l<s, 0<k<s$:
    $q((l,k),(l+1,k-1))=1,~q((l,k),(l+1,k+1))=1,~
    q((l,k),(l+1,k-1))=2$ and $q((l,k),(l-1,k+1))=2$
    \item $l=s, 0<k<s$: $q((l,k),(l-1,k-1))=1$ and
    $q((l,k),(l-1,k+1))=1$
    \item and $0$ otherwise.
\end{itemize}

\begin{figure}[h]

\setlength{\unitlength}{1cm}
\begin{picture}(12,12)

\put(1,1){\vector(1,0){11}}
\put(12,0.5){$k$}

\put(2,1){\vector(0,1){11}}

\put(1.5,0.5){$(1,0)$}

\put(1.5,12){$l$}


\multiput(3,1)(2,0){4}{\circle*{0.1}} \multiput(2,2)(2,0){5}{\circle*{0.1}}
\multiput(3,3)(2,0){5}{\circle*{0.1}} \multiput(4,4)(2,0){4}{\circle*{0.1}}
\multiput(5,5)(2,0){4}{\circle*{0.1}} \multiput(6,6)(2,0){3}{\circle*{0.1}}
\multiput(7,7)(2,0){3}{\circle*{0.1}} \multiput(8,8)(2,0){2}{\circle*{0.1}}
\multiput(9,9)(2,0){2}{\circle*{0.1}} \multiput(10,10)(2,0){1}{\circle*{0.1}}
\multiput(11,11)(2,0){1}{\circle*{0.1}}


\put(8,4){\vector(-1,-1){0.9}} \put(7.5,3.7){\makebox(0,0)[b]{$1$}}

\put(8,4){\vector(-1,1){0.9}}\put(7.5,4.7){\makebox(0,0)[b]{$1$}}

\put(8,4){\vector(1,-1){0.9}}\put(8.5,3.7){\makebox(0,0)[b]{$2$}}

\put(8,4){\vector(1,1){0.9}}\put(8.5,4.7){\makebox(0,0)[b]{$2$}}


\put(7,1){\vector(-1,1){0.9}}\put(6.5,1.7){\makebox(0,0)[b]{$2$}}

\put(7,1){\vector(1,1){0.9}}\put(7.5,1.7){\makebox(0,0)[b]{$4$}}


\put(11,7){\vector(-1,-1){0.9}} \put(10.5,6.7){\makebox(0,0)[b]{$1$}}

\put(11,7){\vector(-1,+1){0.9}}\put(10.5,7.7){\makebox(0,0)[b]{$1$}}

\put(11,1){\vector(-1,+1){0.9}}\put(10.5,1.7){\makebox(0,0)[b]{$2$}}


\put(2,2){\vector(1,1){0.9}}\put(2.5,2.7){\makebox(0,0)[b]{$3$}}

\put(2,2){\vector(1,-1){0.9}}\put(2.5,1.7){\makebox(0,0)[b]{$1$}}


\put(7,7){\vector(1,1){0.9}}\put(7.5,7.7){\makebox(0,0)[b]{$3$}}

\put(7,7){\vector(1,-1){0.9}}\put(7.5,6.7){\makebox(0,0)[b]{$1$}}

\put(7,7){\vector(-1,-1){0.9}}\put(6.5,6.7){\makebox(0,0)[b]{$2$}}


\put(11,11){\vector(-1,-1){0.9}}\put(10.5,10.7){\makebox(0,0)[b]{$2$}}

\end{picture}
\caption{Local configurations and some transition rates for $s=10$.}
\label{figure1}
\end{figure}

The position of the spider will be defined as follows. Assume the two legs be in positions $x_1$ and $x_2$ with
$d(x_1, x_2)=l$. Let $\langle x_1=y_1,y_2,\ldots,y_l=x_2\rangle$  be the geodesic between $x_1$ and $x_2$. If
$l$ is even the position of the spider is defined as $y_{l/2}$ and if $l$ is odd as the middle of the edge
$(y_{(l-1)/2},y_{(l+1)/2})$. Therefore, the distance to the root is $d(y_{l/2},\orig)$ resp.
$d(y_{(l-1)/2},\orig)+1/2$. In order to estimate the speed observe that only local configurations of the type
$(l,0),~1\leq l\leq s$ have a positive drift, i.e., $\E [ H(S_1)-H(S_0)| \tau(S_0)=O_i]>0$ only if $O_i$
corresponds to $(l,0),~1\leq l\leq s$, compare with Example \ref{ex:tree_rooted_at_inf}. Let $\Pi$ the
stationary distribution of the factor chain $(Y_n)_{n\geq 0}$. In order to prove that $V(s)\to 0$ as $s\to \infty$ it
remains to show that $\Pi(\{(l,0), ~1\leq l\leq s\})\to 0$ as $s\to \infty$. Let $\tau_{B}:=\min\{n\geq 1:~
Y_n\in B\}$ be the first hitting time of $B$ and $m_{x,B}=\E_x \tau_{B}$. There is the following relation
\begin{equation*}\label{eq:pi_mib}
\sum_{i\in B} \Pi(i) m_{i,B}=1,
\end{equation*}
see Proposition 6.24 of \cite{KSK}. Let $B=\{(l,0), ~1\leq l\leq s\}$ and observe that the projection of $(Y_n)_{n\geq 0}$
on the second coordinate $k$ behaves like a reflected simple random walk and hence $m_{x,B}\geq c s$ for all $s$
and some constant $c>0.$ Eventually, $\Pi(B)\to 0$ as $s\to\infty$  and the claim follows.
\end{proof}

\subsection*{Some open questions}\label{subs:out}

\begin{que}
Consider a transient SRW on a graph $G$ with positive speed $V$ and denote $V(s)$ the speed of a $k$-leg spider
with span $s$. Is it true that $V(s)\to c< V$ as $s\to \infty$? For which graphs $G$ does $V(s)\to 0$?
\end{que}

\begin{que}
Consider a SRW on a transitive graph with positive speed. Is it true that every transitive spider has positive speed?
\end{que}

\begin{conj}
If the SRW on a Cayley  graph $G$ has positive speed $V$, then any transitive spider has positive speed smaller than $V$.
\end{conj}

\section*{Acknowledgements}
We warmly thank the anonymous referee for his careful reading and numerous suggestions.
C.G. is grateful to Fapesp (grant  2009/51139--3) for financial support. S.M. thanks DFG (project MU 2868/1--1) and Fapesp (grant 2009/08665--6) for financial support. S.P. is grateful to Fapesp (thematic grant 2009/52379--8), CNPq (grants 300886/2008--0,  472431/2009--9) for financial support.  S.M. and S.P. thank CAPES/DAAD (Probral) for support.

\begin{small}
\addcontentsline{toc}{chapter}{Bibliography}
\bibliography{bib_spiders}
\end{small}

\end{document}